\documentclass[12pt,reqno]{amsart}
\usepackage[utf8x]{inputenc}

\usepackage{hyperref}
\hypersetup{
    colorlinks=true,       
    linkcolor=blue,          
    citecolor=blue,        
    filecolor=blue,      
    urlcolor=cyan  
}
\usepackage{latexsym,amssymb,amsmath,amsthm,bbm}
\usepackage{epsfig,mathrsfs}
\usepackage{pst-eucl}
\usepackage{pst-plot,pstricks-add}

\newcommand{\R}{{\ensuremath{\mathbb{R}}}}
\newcommand{\N}{{\ensuremath{\mathbb{N}}}}

\newcommand{\B}{{\ensuremath{\mathcal{B}}}}
\renewcommand{\P}{\ensuremath{\mathbb{P}}}
\renewcommand{\dj}{d\kern-0.4em\char"16\kern-0.1em}
\newcommand{\E}{\ensuremath{\mathbb{E}}}

\newtheorem{Thm}{Theorem}[section]
\newtheorem{Cor}[Thm]{Corollary}
\newtheorem{Lem}[Thm]{Lemma}
\newtheorem{Prop}[Thm]{Proposition}

\theoremstyle{remark}
\newtheorem{Rem}[Thm]{Remark}

\theoremstyle{definition}
\newtheorem{Ex}[Thm]{Example}

\theoremstyle{definition}

\theoremstyle{definition}

\newcommand{\cal}[1]{\mathcal{#1}}

  \def\sL {{\cal L}}

  \def\sU {{\cal U}}

  \def\bR {{\mathbb R}}

\def\wt{\widetilde}

\begin{document}
\numberwithin{equation}{section}
\bibliographystyle{amsalpha}

\title[Harnack inequalities for SBM]{Harnack inequalities for subordinate
Brownian
motions}
\begin{abstract}
In this paper, we consider 
subordinate
 Brownian motion 
 $X$ in $\R^d$, $d \ge 1$, 
 where the Laplace exponent $\phi$
of the corresponding subordinator satisfies some mild conditions. 
The scale
invariant Harnack inequality 
 is proved for $X$. 
  We first give new forms of asymptotical
properties of the L\' evy and potential density of the subordinator near zero. 
Using these results 
we find asymptotics of the L\' evy density
and potential density of $X$ near the origin, which is essential to our approach. The examples which are covered by
our results
include geometric stable processes and relativistic geometric stable processes,
i.e. the cases when the subordinator has the Laplace exponent
\[
\phi(\lambda)=\log(1+\lambda^{\alpha/2})\ \ \ \ 
(0<\alpha\leq 2)\]
and
\[
\phi(\lambda)=\log(1+(\lambda+m^{\alpha/2})^{2/\alpha}-m)\ \ \ \ (0<\alpha<2,\,
m>0)\,.
\]
\end{abstract}

\author{Panki Kim}
\address{Department of Mathematical Sciences and Research Institute of Mathematics,
Seoul National University,
Building 27, 1 Gwanak-ro, Gwanak-gu
Seoul 151-747, Republic of Korea}
\curraddr{}
\email{pkim@snu.ac.kr}
\thanks{Research supported by Basic Science Research Program through
the National Research Foundation of Korea(NRF) funded by the Ministry of Education,
Science and Technology(0409-20110087).}

\author{Ante Mimica}
\address{Fakult\"{a}t f\"{u}r Mathematik, Universit\"{a}t Bielefeld, Postfach
100131, D-33501 Bielefeld, Germany}
\curraddr{}
\thanks{ Research supported in part by German Science Foundation DFG via IGK
"Stochastics and real world models" and SFB 701.}
\email{amimica@math.uni-bielefeld.de}
\thanks{}

\subjclass[2010]
{Primary 60J45, Secondary 60G50, 60G51, 60J25, 60J27}

\keywords{geometric stable process, Green function, Harnack inequality, Poisson
kernel, harmonic function, potential, subordinator, subordinate
Brownian motion}

\maketitle

\section*{Introduction}
Consider a Brownian motion  $B=(B_t,\P_x)$ in $\R^d$, $d \ge 1$,  and an independent
subordinator $S=(S_t\colon t\geq 0)$. It is known that the stochastic
process $X=(X_t,\P_x)$ defined by $X_t=B(S_t)$ is a L\' evy process. The
process $X$ is called the subordinate Brownian motion. 

A  non-negative function $h\colon \R^d\rightarrow
[0,\infty)$ is said to be harmonic 
with respect to $X$ in an open set
$D\subset\R^d$ if for all open sets $B\subset \R^d$ 
whose closure is compact and  
contained in
$D$ the following mean value property holds
\[
h(x)=\E_x[h(X_{\tau_B})]\ \ \text{ for all } \ \ x\in B\,,
\]
where $\tau_B=\inf\{t>0\colon X_t\not\in B\}$ denotes the first exit time from
the set $B$ .

The Harnack inequality holds for the process $X$ if there exists a constant
$c>0$ such that for any $r\in (0,1)$ and any non-negative function $h$ on $\R^d$
which is harmonic in ball $B(0,r)=\{z\in \R^d\colon |z|<r\}$ the following
inequality is true
\begin{equation}\label{eq:intro-1}
h(x)\leq c\, h(y)\ \ \text{ for all }\ \ x,y\in B(0,\tfrac{r}{2})\,.
\end{equation}

Space homogeneity of L\' evy processes implies that the same inequality is true
on any ball $B(x_0,r)=\{z\in
\R^d\colon |z-x_0|<r\}$. This type of Harnack inequality is sometimes called
scale invariant (or geometric) Harnack inequality, since the constant $c$ in
(\ref{eq:intro-1}) stays the same for any $r\in (0,1)$.

The main goal of this paper is to prove the scale invariant Harnack inequality
for a class of subordinate Brownian motions.
Our most important
contribution is that within our framework we can treat subordinate
Brownian motions with subordinators whose Laplace exponent
\[\phi(\lambda):=-\log \E e^{-\lambda S_t}\] varies slowly at infinity. In
particular, we are able to give a positive answer for many processes for which
only the non-scale version of Harnack inequality was known so far.

Here are few examples of such processes.

{\bf Example 1} (Geometric stable processes)
\[\phi(\lambda)=\log(1+\lambda^{\beta/2}),\ \  \ (0<\beta\leq 2 
).\]

{\bf Example 2} (Iterated geometric stable processes)
\begin{align*}
 \phi_1(\lambda)&=\log(1+\lambda^{\beta/2})  \ \ \ (0<\beta\leq 2)\\
\phi_{n+1}&=\phi_1\circ \phi_n \ \ \  n\in \N.
\end{align*}

{\bf Example 3} (Relativistic geometric stable processes)
\[\phi(\lambda)=\log\left(1+\left(\lambda+m^{\beta/2}\right)^{2/\beta}-m
\right)  \ \ \ (m>0,\,0<\beta<2
).\]

\begin{Rem}
 The non-scale version of Harnack inequality for
geometric stable and iterated geometric stable processes was proved in
\cite{SSV2}. It was not known
whether scale invariant version of this inequality held. Recently this turned
out to
be the case in dimension $d=1$ 
(see \cite{GR2}). In \cite{GR2} the authors 
used theory of fluctuation of one-dimensional L\'
evy processes and it was not clear how to generalize this technique to higher
dimensions. 
Nevertheles, this result suggests that the scale invariant version of Harnack
inequality may hold in higher dimensions.
\end{Rem}

Another feature of our approach is that it is unifying in the
sense that it covers many classes of subordinate Brownian motions for which
the scale invariant Harnack inequality was recently proved. For example we can
treat many subordinators whose Laplace exponent varies regularly at infinity.
As a special example,  rotationally
invariant $\alpha$-stable processes ($\alpha\in (0,2)$) are included in  our
framework.

Let us be more precise now.  In this paper we consider subordinate Brownian
motions 
$X$ in $\R^d$ $(d \ge 1$), for which the
Laplace exponent $\phi$ of the corresponding subordinator $S$ satisfies (see
Sections \ref{sec:prelim} and \ref{sec:2}
 for details concerning these conditions):
\begin{itemize}
	\item[{\bf (A-1)}] 
	the potential measure of $S$ has a decreasing density;
	\item[{\bf(A-2)}] 
	the L\' evy measure of $S$ is infinite and has a decreasing density;
	\item[{\bf(A-3)}] there exist constants $\sigma>0$, $\lambda_0>0$ and
$\delta \in (0, 1]$
 such that 
\begin{equation*}
  \frac{\phi'(\lambda x)}{\phi'(\lambda)}\leq
\sigma\,x^{-\delta}\ \text{ for all
}\ x\geq 1\ \text{ and }\ \lambda\geq\lambda_0\,;
\end{equation*}
\end{itemize}

Our main result is the following scale invariant Harnack inequality. 
\begin{Thm}[Harnack inequality]\label{thm:harnack}
Suppose $X$ is a subordinate Brownian motion satisfying {\bf (A-1)}--{\bf
(A-3)}. We further assume that 
the L\' evy density $J(x)=j(|x|)$ of $X$ 
satisfies 
\begin{equation}\label{eq:sub-11}
 j(r+1)\leq j(r)\leq c'j(r+1),\ 
r>1,
\end{equation} 
for some constant $c'\geq 1$.

	Then	there exists a constant $c>0$  such that for
all $x_0\in \R^d$ and $r\in (0,1)$
	\[
		h(x_1)\leq c\, h(x_2)\ \text{ for all }\ x_1,x_2\in
B(x_0,\tfrac{r}{2})
	\]
	and for every  non-negative function $h\colon \R^d\rightarrow
[0,\infty)$ which
is harmonic in $B(x_0,r)$.
\end{Thm}

As already mentioned at the beginning, this theorem is a new result for
{ Examples 1--3} above.  Note that  when $\phi$ is a complete Bernstein
function, under the assumption {\bf (A-2)},  {\bf (A-1)} and \eqref{eq:sub-11}
hold (see Corollary
10.7 in \cite{SSV} and our Remark \ref{r:sd}).

The condition {\bf (A-3)} is implied by the following stronger condition
  \begin{equation}\label{eq:stri-cond}
    \forall\, x>0\ \ \ \ \lim_{\lambda\to\infty}\frac{\phi'(\lambda
x)}{\phi'(\lambda)}=x^{\frac{\alpha}{2}-1}\ \ \ \ (0\leq
\alpha<2)\,.
  \end{equation}
In other words, (\ref{eq:stri-cond}) says that
$\phi'$ varies regularly at infinity with index $\frac{\alpha}{2}-1$.
{\ Examples  1--3} satisfy this condition with $\alpha=0$.

As already mentioned, Theorem \ref{thm:harnack} covers also processes for which
the Harnack inequality was known before
(see \cite{KSV3,Mi, RSV}):

{\bf Example 4}\ \ Assume that $\phi$ satisfies {\bf (A-1)}, {\bf (A-2)} and
\[
 \phi(\lambda)\asymp \lambda^{\alpha/2}\ell(\lambda),\ \lambda \to\infty\ \ \
(0<\alpha<2
)\, 
\]
where $\ell$ varies slowly at infinity, i.e.
\[
 \forall\, x>0\ \ \lim_{\lambda \to\infty}\frac{\ell(\lambda
x)}{\ell(\lambda)}=1
\]
($f(\lambda) \asymp g(\lambda),\ \lambda \to \infty$ means that
$f(\lambda)/g(\lambda)$ stays bounded between two positive constants as $\lambda
\to \infty$ ).
We can take, for example,  $\ell(\lambda)=[\log(1+\lambda)]^{1-\alpha/2}$ or
$\ell(\lambda)=[\log(1+\log(1+\lambda))]^{1-\alpha/2}$.

The main ingredient in our proof of Harnack inequality is a good
estimate of the Green function $G_{B(0,r)}(x,y)$ of the ball $B(0,r)$ when $y$
is near its boundary. To be more precise, we will prove that there are 
a function $\xi\colon (0,1)\rightarrow (0,\infty)$
and constants
 $c_1,c_2>0$ and $0<\kappa_1<\kappa_2<1$ such that for every $r\in
(0,1)$, 
\begin{equation}\label{eq:intro-xi}
 c_1 \xi(r)r^{-d}\,\E_y\tau_{B(0,r)}\leq G_{B(0,r)}(x,y)\leq c_2
\xi(r)r^{-d}\,\E_y\tau_{B(0,r)},
\end{equation}
for $x\in B(0,\kappa_1 r)$ and $y\in B(0,r)\setminus B(0,\kappa_2r)$ (see
Corollary \ref{cor:green})\,. 

Depending on the considered process, the function 
$r \mapsto \xi(r)$ 
can have two different types of
behavior. 
For example, it turns out
that in { Example 1}
  \[\xi(r)\asymp \tfrac{1}{\log(r^{-1})}\ \text{ as }\ r\to
0+\,,\]
while in { Example 4}
\[
 \xi(r)\asymp1\ \text{ as }\ r\to
0+\,.
\]

To obtain the mentioned estimates of the  Green function we use
asymptotical properties of L\' evy density $\mu(t)$ and potential density $u(t)$
of the underlying subordinator near zero. It turns out that it is not possible
 use Tauberian theorems in each case. In  Section \ref{sec:2} we obtain needed
asymptotical properties without use of such theorems. The asympotical behavior
can be expressed in the terms of the Laplace exponent in the following way
\[
 \mu(t)\asymp t^{-2}\phi'(t^{-1})\ \ \text{ and }\ \ \ u(t)\asymp
t^{-2}\frac{\phi'(t^{-2})}{\phi(t^{-2})^2}\ \ \text{ as }\ \ t\to 0+\,.
\]

Harnack inequalities for symmetric stable L\' evy processes were
obtained in
\cite{BSS, BS}.
 A new technique on Harnack inequalities for stable like jump processes
was developed in \cite{BL1} and generalized in \cite{SV1}. Similar technique
was used for various jump processes in 
\cite{BKa, CK2,CK}.
 In \cite{KS} the
Harnack inequality was proved for truncated stable processes and it was
generalized in \cite{Mi}. Harnack inequality for some classes subordinate
Brownian motions was also considered in \cite{KSV3}.

Let us comment what happens when 
one applies  techniques
developed for jump processes (as in \cite{BL1}) to our situation.
The proof in this case relied on an
estimate of Krylov-Safonov type: 
there exists a constant $c>0$ such
that 
\[
	\P_x(T_A<\tau_{B(0,r)})\geq c\,\frac{|A|}{|B(0,r)|}
\]
for any $r\in (0,1)$, $x\in B(0,\frac{r}{2})$ and $A\subset \R^d$ closed, where
$T_A=\tau_{A^c}$ denotes the first hitting time of the set $A$ and $|A|$
denotes its Lebesgue
measure.

Although this technique is quite general and can be applied to a  much larger
class of Markov jump processes, there are some situations in our setting which
show that it is not
 applicable even to a rotationally invariant L\' evy process.
A good example is the proof  of the Harnack
inequality in \cite{SSV2}, where the mentioned Krylov-Safonov type estimate was
indispensable. Contrary to the case of stable-like
processes, this
estimate is not uniform in $r\in (0,1)$. 

For example, for a geometric stable
process it is possible to find a sequence of radii $(r_n)$ and closed sets
$A_n \subset B(0,r_n)$ such that $r_n\to 0$, $\frac{|A_n|}{|B(0,r_n)|}\geq
\frac{1}{4} $ and 
\[
  \P_0(T_{A_n}<\tau_{B(0,r_n)})\to 0,\ \text{ as }\ n\to\infty\,. 
\]
This non-uniformity does not allow to obtain scale-invariant Harnack inequality
using this technique. 
In this sense, the investigation of Harnack inequality becomes interesting even
in
the case of a L\' evy process. We have not encountered a technique so far that
would cover cases of a more general jump process in this direction.

The paper is organized as follows. In Section \ref{sec:prelim} we give basic
notions which we use in sections that follow. A new forms of asymptotical
properties
of the L\' evy and the potential densities of subordinators are obtained in
Section \ref{sec:2}. Technical lemmas concerning 
asymptotic inversion of the Laplace transform used in this section are deferred
to Appendix \ref{app:a}. These results in Appendix \ref{app:a} can be also of
independent interest, since they represent an alternative to the Tauberian
theorems, which were mainly used in previous 
works.

Using results of the Section \ref{sec:2} we obtain the behavior of the L\' evy
measure and the Green function (potential) of the process $X$ in Section
\ref{sec:3}. In Section \ref{sec:4} we obtain pointwise estimates of the Green
functions of small balls needed to prove the main result, which is proved in
Section \ref{sec:5}.

{\bf Notation. } Throughout the paper we use the notation $f(r)\asymp
g(r),\ r\to a$ to denote that ${f(r)}/{g(r)}$ stays between two positive
constants as $r\to a$. 
Simply, $f\asymp g$  means that the quotient $f(r)/g(r)$ stays bounded
between two positive numbers on their common domain of definition. 
We say that $f\colon\R\rightarrow \R$ is increasing if
$s\leq t$ implies $f(s)\leq f(t)$ and analogously for a decreasing function.  
For a Borel
set $A\subset \bR^d$, we also use $|A|$ to denote its Lebesgue
measure. We will use ``$:=$" to denote a definition,
which is read as ``is defined to be".
 For any $a, b\in \bR$, we use the notations
$a\wedge b:=\min \{a, b\}$ and $a\vee b:=\max\{a, b\}$.
 The values of the constants $c_1, c_2, \cdots$ stand for constants whose values
are unimportant and which may change from location to location.
The labeling of the constants $c_1, c_2, \ldots$ starts anew in the proof of
each result. 

\section{Preliminaries}\label{sec:prelim}

A stochastic process  $X=(X_t,\P_x)$ in $\R^d$ is said to be a pure jump L\'
evy 
process if it has stationary and independent increments, its trajectories
are 
right-continuous with left limits and the 
characteristic exponent $\Phi$ in 
\[
\E_x\left[\exp{\{i\langle \xi,X_t-X_0\rangle\}}\right]=\exp{\{-t\Phi(\xi)\}},\
\xi\in \R^d
\]
is of the form
\begin{equation}\label{eq:sub-1}
\Phi(\xi)=\int_{\R^d}\left(1-\exp{\{i\langle\xi,
x\rangle\}}+i\langle\xi, x\rangle\, 1_{\{|x|<1\}}\right)\Pi(dx).
\end{equation}
The measure $\Pi$ in (\ref{eq:sub-1}) is called the L\' evy measure of $X$ and
it satisfies
$
 \Pi(\{0\})=0$  and $ \int_{\R^d}(1\wedge |x|^2)\Pi(dx)<\infty\,
.
$

Let  $S=(S_t\colon t\geq 0)$ be a subordinator, i.e. a 
L\' evy process taking values in $[0,\infty)$ and starting at $0$. It is
more convenient to consider the Laplace transform in this case
\begin{equation}\label{eq:sub-2}
\E \exp{\{-\lambda S_t\}}=\exp{\{-t\phi(\lambda)\}}.
\end{equation}
The function $\phi$ in (\ref{eq:sub-2}) is called the Laplace exponent of $S$
and it is of the form
\begin{equation}\label{eq:sub-3}
\phi(\lambda)=\gamma t+\int_{(0,\infty)}(1-e^{-\lambda t})\,\mu(dt)\,,
\end{equation}
where $\gamma\geq 0$ and the L\' evy measure 
$\mu$ of $S$ is now a measure on
$(0,\infty)$
satisfying $\int_{(0,\infty)}(1\wedge t)\mu(dt)<\infty$ (see p. 72 in 
\cite{Be}).

The function $\phi$ is an example of a Bernstein function, i.e.
$\phi\in C^\infty(0,\infty)$ and $(-1)^n\phi^{(n)}\leq 0$ for all
$n\in \N$ (see p. 15 in \cite{SSV}). Here $\phi^{(n)}$ denotes the $n$-th
derivative of $\phi$. Conversely, every Bernstein function $\phi$ satisfying
$\phi(0+)=0$ has a representation (\ref{eq:sub-3}) and there exists a
subordinator with the Laplace exponent $\phi$ .

The potential measure of the subordinator $S$ is defined by
\begin{equation}\label{eq:sub-8}
 U(A)=\int_0^\infty\P(S_t\in A)\,dt.
\end{equation}
The Laplace transform of $U$ is then
\begin{equation}\label{eq:intro-lp}
\sL U(\lambda)=\E\int_{(0,\infty)}e^{-\lambda
S_t}\,dt=\frac{1}{\phi(\lambda)},\, \quad \lambda>0\,.
\end{equation}

A Bernstein function $\phi$
is said to be a complete Bernstein function if the L\' evy measure $\mu$ has a
completely monotone density, i.e. $\mu(dt)=\mu(t)\,dt$ with $\mu\in
C^\infty(0,\infty)$ satisfying $(-1)^n\mu^{(n)}\geq 0$ for all $n\in
\N\cup\{0\}$\,. The corresponding subordinator is called a complete
subordinator. 

In this case we can control large jumps of L\' evy density
$\mu$ in the following way. There exists a constant $c>0$ such that
\begin{equation}\label{eq:sub-4}
\mu(t)\leq c\,\mu(t+1)\ \text{ for all }\ t\geq 1
\end{equation}
(see Lemma 2.1 in \cite{KSV2}).
If, in addition,  $\mu(0,\infty)=\infty$, the potential measure $U$ has a
decreasing density, i.e. there exists a decreasing function $u\colon
(0,\infty)\rightarrow (0,\infty)$ such that $U(dt)=u(t)\,dt$ (see Corollary
10.7 in \cite{SSV}).

Let $B=(B_t,\P_x)$ be a Brownian motion in $\R^d$
(running with a time clock twice as fast as the standard Brownian motion)
 and let $S=(S_t\colon t\geq
0)$ be an independent subordinator. 
We define a new process $X=(X_t,\P_x)$ by $X_t=B(S_t)$ and call it subordinate
Brownian motion. This process is a L\' evy process with the characteristic
exponent $\Phi(\xi)=\phi(|\xi|^2)$. Moreover, $\Phi$ has representation
(\ref{eq:sub-1}), with the L\' evy measure of the form $\Pi(dx)=j(|x|)\,dx$
and
\begin{equation}\label{eq:sub-5}
 j(r)=\int_{(0,\infty)} (4\pi
t)^{-d/2}\exp{
\left(-\tfrac{r^2}{4t}\right)}\mu(dt),\, r>0
\end{equation}
(see Theorem 30.1 in \cite{S}).

The process $X$ has a transition density $p(t,x,y)$ given by
\begin{equation}\label{eq:sub-6}
 p(t,x,y)=\int_0^\infty (4\pi
t)^{-d/2}\exp{\left(-\tfrac{|x-y|^2}{4t}\right)}\P(S_t\in ds)\,.
\end{equation}

The process $X$ is said to be transient if
$\P_0(\lim_{t\to\infty}|X_t|=\infty)=1$. 
Since the characteristic exponent of
$X$ is symmetric we have the following Chung-Fuchs type criterion for transience
\begin{align}\nonumber
	X\ \ \text{ is transient } \ \ &\iff\ \ \
\int_{B(0,R)}\frac{d\xi}{\phi(|\xi|^2)}<\infty \ \text{ for some }\ R>0\\
	&\iff\ \ \
\int_0^R\frac{\lambda^{\frac{d}{2}-1}}{\phi(\lambda)}\,d\lambda<\infty \ \text{
for some }\ R>0\label{eq:prelim-10}\\
	&\iff \E_0\left[\int_0^\infty 1_{\{|X_t|<R\}}\,dt\right]<\infty \ \
\text{ for every }\ R>0\label{eq:prelim-255}
\end{align}
(see Corollary 37.6 and Theorem 35.4 in \cite{S}). 

In this case  we can define the Green function (potential) by
$
G(x,y)=\int_0^\infty p(t,x,y)\,dt\,.
$
Then \eqref{eq:sub-8} and \eqref{eq:sub-6} 
give us a useful formula
$G(x,y)=G(y-x)=g(|y-x|)$, where
\begin{equation}\label{eq:sub-7}
 g(r)=\int_{(0,\infty)} (4\pi
t)^{-d/2}\exp{\left(-\tfrac{r^2}{4t}\right)}U(dt),\ r>0\,.
\end{equation}
Note that $g$ and $j$ are decreasing.

Let $D\subset \R^d$ be a bounded open subset. We define killed process $X^D$ by
$X^D_t=X_t$ if $t<\tau_D$ and $X^D_t=\Delta$ otherwise, where $\Delta$ is some
point adjoined to $D$ (usually called cemetery).

The transition density and the Green function of $X^D$ are given by
\[
 p_D(t,x,y)=p(t,x,y)-\E_x\left[p(t-\tau_D,X(\tau_D),y); \tau_D<t\right]
\]
and $G_D(x,y)=\int_0^\infty p_D(t,x,y)\,dt$. In the transient case we have the
following formula
\begin{equation}\label{eq:sub-9}
 G_D(x,y)=G(x,y)-\E_x[G(X(\tau_D),y)]\,.
\end{equation}
Also, $G_D(x,y)$ is symmetric and, for fixed $y\in D$, $G_D(\cdot,y)$ is 
harmonic in $D\setminus\{y\}$ .
Furthermore,  $G_D\colon (D\times D)\setminus \{(x,x)\colon x\in D
\}\rightarrow [0,\infty)$ and $x\mapsto \E_x\tau_D$ are continuous functions. 

By the result of Ikeda and Watanabe (see Theorem 1 in \cite{IW})
\begin{equation}\label{eq:sub-105}
 \P_x(X_{\tau_D}\in F)=\int_F\int_DG_D(x,y)j(|z-y|)\,dy\,dz
\end{equation}
for any $F\subset \overline{D}^c$. 
If we define the Poisson kernel of the set $D$ by
\begin{equation}\label{eq:sub-10} 
 K_D(x,z)=\int_D G_D(x,y)j(|z-y|)\,dy,
\end{equation}
then $\P_x(X_{\tau_D}\in F)=\int_F K_D(x,z)\,dz$ for any $F\subset
\overline{D}^c$. In other words, the Poisson kernel is the density of the exit
distribution. 

Since a subordinate Brownian motion is a rotationally invariant L\' evy process,
it follows
that in the case of the subordinator with zero drift
\[
	\P_x(X_{\tau_{B(x_0,r)}}\in \partial B(x_0,r))=0
\]
(see \cite{Sz}) 
and thus, 
for a measurable function  $h\colon \R^d\rightarrow [0,\infty)$ ,
\begin{equation}\label{e:new1}
	\E_x[h(X_{\tau_{B(z_0,s)}})]=\int_{\overline{B(z_0,s)}^c}K_{B(z_0,s)}(x,z)h(z)\,dz
\end{equation}
for any ball $B(z_0,s)$ .

\section{Subordinators}\label{sec:2}

Let $S=(S_t\colon t\geq 0)$ be a subordinator with the Laplace exponent $\phi$
satisfying the following conditions:
\begin{itemize}
	\item[{\bf (A-1)}] 
	the potential measure $U$ of $S$ has a decreasing density $u$. i.e., there is a decreasing function
$u\colon (0,\infty)\rightarrow (0,\infty)$ so that $U(dt)=u(t)\,dt$;
	\item[{\bf(A-2)}] 
the L\' evy measure $\mu$ of $S$ is infinite, i.e. $\mu(0,\infty)=\infty$,  and has a decreasing density $t \to \mu(t)$;
	\item[{\bf(A-3)}] there exist constants $\sigma>0$, $\lambda_0>0$ and
$\delta>0$ such that 
\begin{equation}\label{e:new21}
  \frac{\phi'(\lambda x)}{\phi'(\lambda)}\leq
\sigma\,x^{-\delta}\ \text{ for all
}\ x\geq 1\ \text{ and }\ \lambda\geq\lambda_0\,
\end{equation}
\end{itemize}

\begin{Rem}
	\begin{itemize}
		\item[(i)] {\bf (A-1)} and {\bf (A-2)} imply that 
$\phi$ is a special Bernstein function, 
i.e., $\lambda \to \lambda/\phi(\lambda)$ is also a Bernstein function (see pp. 92-93 in
\cite{SSV}).
		\item[(ii)] 
	{\bf(A-2)} implies that $$
			\phi(\lambda)=\gamma \lambda+\int_0^\infty (1-e^{-\lambda
t})\mu(t)dt.$$
		Thus {\bf(A-2)}--{\bf (A-3)} imply $\gamma=0$, by letting $x\to +\infty$ .
	\end{itemize}
\end{Rem}

First we prove a simple result that holds for any Bernstein function, which 
will be used in Section \ref{sec:4}.
\begin{Lem}\label{lem:sub-50} Let $\phi$ be a Bernstein function. 
\begin{itemize}
\item[(i)]
For every $x\geq 1$, 
$
 \phi(\lambda x)\leq x \phi(\lambda)\ \text{ for all }\ 
 \lambda >0\,.
$
\item[(ii)] 
\eqref{e:new21} implies that 
for every $\varepsilon>0$ there is a constant $c=c(\varepsilon)>0$
so that
\[
	\frac{\phi(\lambda x)}{\phi(\lambda)}\leq c\,x^{1-\delta +\varepsilon}\
\text{ for all }\ \lambda\geq \lambda_0 \ \text{ and }\ x\geq 1\,.
\]
\end{itemize}
\end{Lem}
\proof (i)
Since $\phi'$ is decreasing and $x\geq 1$, 
\begin{align*}
 \phi(\lambda x)=\int_0^{\lambda x}\phi'(s)\,ds \leq
\int_0^{\lambda x}\phi'(\tfrac{s }{x})\,ds=x\phi(\lambda )\,.
\end{align*}

(ii) Without loss of generality we may assume that $\sigma\geq 2$ in {\bf (A-3)}. Using {\bf (A-3)},
for any $k\geq 2$ the following recursive inequality holds
\begin{align*}
	\phi(\lambda \sigma^{\frac{k}{\varepsilon}})-\phi(\lambda
\sigma^{\frac{k-1}{\varepsilon}})&=\int_{\lambda
\sigma^{\frac{k-1}{\varepsilon}}}^{\lambda
\sigma^{\frac{k}{\varepsilon}}}\phi'(s)\,ds
	\leq \sigma^{1-\frac{\delta}{\varepsilon}}\int_{\lambda
\sigma^{\frac{k-1}{\varepsilon}}}^{\lambda
\sigma^{\frac{k}{\varepsilon}}}\phi'(s\sigma^{-\frac{1}{\varepsilon}})\,ds\\
	&=\sigma^{1+\frac{1-\delta}{\varepsilon}}\left(\phi(\lambda
\sigma^{\frac{k-1}{\varepsilon}})-\phi(\lambda
\sigma^{\frac{k-2}{\varepsilon}})\right)\,.
\end{align*}
Iteration yields
\begin{equation}\label{eq:sbm-iter}
	\phi(\lambda \sigma^{\frac{k}{\varepsilon}})-\phi(\lambda
\sigma^{\frac{k-1}{\varepsilon}})\leq \sigma^{(k-1)
(1+\frac{1-\delta}{\varepsilon})}\left(\phi(\lambda
\sigma^{\frac{1}{\varepsilon}})-\phi(\lambda)\right) \text{ for every }k\geq 2.
\end{equation}

Let $n\in \N$ be chosen so that $\sigma^{\frac{n-1}{\varepsilon}}\leq
x<\sigma^{\frac{n}{\varepsilon}}$ .

If $n=1$, then by (i),
$
	\phi(\lambda\sigma^{-\frac{1}{\varepsilon}})\leq
\sigma^{\frac{1}{\varepsilon}}\phi(\lambda)\leq
\sigma^{\frac{1}{\varepsilon}+\frac{2\delta}{\varepsilon}}x^{-\delta}
\phi(\lambda)
$
which, by monotonicity of $\phi$, implies that $\frac{\phi(\lambda
x)}{\phi(\lambda)}\leq \sigma^{\frac{1+2\delta}{\varepsilon}}x^{-\delta}$ .

Let us consider now the case $n\geq 2$. Using (\ref{eq:sbm-iter}) and (i) we
deduce
\begin{align*}
	&\phi(\lambda
\sigma^{\frac{n}{\varepsilon}})-\phi(\lambda)=\left(\phi(\lambda
\sigma^{\frac{1}{\varepsilon}})-\phi(\lambda)\right)\, \sum_{k=2}^n\sigma^{(k-1)
(1+\frac{1-\delta}{\varepsilon})}\\
	&\leq \left(\phi(\lambda
\sigma^{\frac{1}{\varepsilon}})-\phi(\lambda)\right)\frac{\sigma^{n
(1+\frac{1-\delta}{\varepsilon})}}{\sigma^{1+\frac{1-\delta}{\varepsilon}}-1} \leq \sigma^{\frac{1}{\varepsilon}}\phi(\lambda)\sigma^{n
(1+\frac{1-\delta}{\varepsilon})}\,.
\end{align*}

Therefore
\begin{align*}
	\phi(\lambda x)&\leq \phi(\lambda\sigma^{\frac{n}{\varepsilon}})\leq
2\sigma^{1+\frac{2-\delta}{\varepsilon}}\phi(\lambda)\left(\sigma^{\frac{n-1}{\varepsilon
}}\right)^{\varepsilon +1-\delta}\leq
2\sigma^{1+\frac{2-\delta}{\varepsilon}}\phi(\lambda)x^{\varepsilon+1-\delta}\,.
\end{align*}
\qed

\begin{Prop}\label{prop:sub-levy}
 If {\bf (A-2)}
and {\bf (A-3)} hold, then
\[
\mu(t)\asymp t^{-2}\phi'(t^{-1}),\ t\to 0+\,.
\]
\end{Prop}
\proof
Note that  
\[
 \phi(\lambda+\varepsilon)-\phi(\varepsilon)=\int_0^\infty
\left(e^{-\lambda t}-e^{-\lambda (t+\varepsilon)}\right)\mu(t)\,dt
\]
for any $\lambda>0$ and $\varepsilon>0$ and thus  the condition
(\ref{eq:cond_main}) in Appendix \ref{app:a} holds with $f=\phi$ and $\nu=\mu$.
Since
$\phi$ is a
Bernstein function, it follows that $\phi'\geq 0$ and $\phi'$ is decreasing.
Now we can apply 
Lemmas \ref{lem:1} and \ref{lem:2}.
\qed

\begin{Prop}\label{prop:sub-pot}
If {\bf (A-1)}
and {\bf (A-3)} hold, then
\[
		u(t)\asymp t^{-2}\frac{\phi'(t^{-1})}{\phi(t^{-1})^2},\ t\to
0+\,.
	\]
\end{Prop}
\proof
	By  (\ref{eq:intro-lp}), with $\psi(\lambda)=\frac{1}{\phi(\lambda)}$, we have 
	$
		\int_0^\infty e^{-\lambda t} u(t)\,dt=\psi(\lambda).
	$
	 Note that, for $\lambda
\geq
\lambda_0$ and $x\geq 1$,  {\bf (A-3)} implies
	\[
		\left|\frac{\psi'(\lambda
x)}{\psi'(\lambda)}\right|=\left(\frac{\phi(\lambda)}{\phi(\lambda
x)}\right)^2\,\frac{\phi'(\lambda x)}{\phi'(\lambda)}\leq \frac{\phi'(\lambda
x)}{\phi'(\lambda)}\leq cx^{-\delta},
	\]
	since $\phi$ is increasing. 
	
	We see that (\ref{eq:cond_main}) in Appendix \ref{app:a} is satisfied
with
$f=\tfrac{1}{\phi}$
and $\nu=u$. Since $\phi$ is a Bernstein function, $\phi'\geq 0$ and $\phi'$ is
a decreasing function. Thus $|f'|=\frac{\phi'}{\phi^2}$ is also a
decreasing
function. The result follows now from 
Lemmas \ref{lem:1} and \ref{lem:2}.
\qed

\section{L\' evy density and potential}\label{sec:3}

In Section \ref{sec:2} we have established asymptotic behavior of the L\' evy and
potential density of $S$ near zero. In this section we are going to use these
results to 
give  new forms of  
asymptotic behavior of the L\' evy density and
potential of the process $X$ near the origin.
Throughout the remainder of the paper,  $X$ is the subordinate
Brownian motion with the characteristic
exponent $\phi(|\xi|^2)$ where $\phi$ is the Laplace exponent of $S$.  
\begin{Lem}\label{lem:gr-10}
Suppose that $\phi$ is a special Bernstein function, 
i.e., $\lambda \to \lambda/\phi(\lambda)$ is also a Bernstein function.
Then  the  functions $\eta_1,\eta_2\colon (0,\infty)\rightarrow (0,\infty)$ given by
\[
 \eta_1(\lambda)=\lambda^2\phi'(\lambda) \ \text{ and }\
\eta_2(\lambda)=\lambda^2\frac{\phi'(\lambda)}{\phi^2(\lambda)}
\]
are increasing.
\end{Lem}
\proof
It is enough to prove that $\eta_2$ is increasing, because
$\eta_1=\eta_2\cdot \phi^2$ is then a product of two increasing functions.

Since $\phi$ is a special Bernstein function, 
\[
 \frac{\lambda}{\phi(\lambda)}=\theta+\int_{(0,
\infty)}(1-e^{-\lambda t})\nu(dt),
\]
for some $\theta\geq 0$ and a L\' evy measure $\nu$ (see pp. 92-93 in
\cite{SSV}).
Then
\begin{align*}
\lambda^2\frac{\phi'(\lambda)}{\phi(\lambda)^2}&=\lambda\left(-\frac{\lambda}{
\phi(\lambda)}\right)'+\frac{\lambda}{ \phi(\lambda)}\\
&=\theta+\int_{(0,\infty)}\left(1-(1+\lambda t)e^{-\lambda t}\right)\nu(dt)\,.
\end{align*}
Now the claim follows, since $\lambda\mapsto 1-(1+\lambda t)e^{-\lambda t}$ is
increasing for any $t> 0$.
\qed

\begin{Prop}\label{prop:pot-11}
If {\bf (A-2)}
and {\bf (A-3)} hold, then
\[
 j(r)\asymp r^{-d-2}\phi'(r^{-2}),\ r\to 0+\,.
\]
\end{Prop}
\proof
We use formula (\ref{eq:sub-5}), i.e.
\[
 j(r)=\int_0^\infty (4\pi
t)^{-d/2}\exp{\left(-\tfrac{r^2}{4t}\right)}\mu(t)\,dt\,.
\]
Proposition
\ref{prop:sub-levy} implies that $\mu(t)\asymp t^{-2}\phi'(t^{-1}),\ t\to 0+$. 

We are going to use Proposition \ref{prop:gr-1} in Appendix \ref{app:a} 
with $A=2$, 
$\eta=\mu$ and
$\psi=\phi'$. In order to do this, we need to check conditions 
(a), (b) and (c)-(ii). The condition (a)
follows from the fact that $\phi$ is a Bernstein function and  Lemma
\ref{lem:gr-10}, while (b) follows from
\[
 \int_1^\infty t^{-d/2}\mu(t)\,dt\leq \int_1^\infty
\mu(t)\,dt=\mu(1,\infty)<\infty,
\]
since $\mu$ is a L\' evy measure.
Finally the condition (c)-(ii) follows from \eqref{e:new22}--\eqref{e:new23}.
\qed

\begin{Rem}\label{r:sd}
If $\phi$ is a complete Bernstein function,  using Lemma 4.2 in
\cite{RSV}, our (\ref{eq:sub-4}) and Proposition \ref{prop:pot-11},  we see that \eqref{eq:sub-11}
 holds.
 \end{Rem}

\begin{Lem}\label{lem:pot-cond}
If {\bf (A-1)}  hold and $X$ is transient, then
	\[
		\int_1^\infty t^{-d/2}u(t)\,dt<\infty\,. 
	\]
\end{Lem}
\proof
It follows from  (\ref{eq:intro-lp}) and (\ref{eq:prelim-10}) that for any $t\geq 1$ and $ R>0$ 
\begin{align*}
\infty
&>\int_0^R\frac{\lambda^{\frac{d}{2}-1}}{\phi(\lambda)}\,
d\lambda=\int_0^R\int_0^\infty\lambda^{\frac{d}{2}-1}e^{-\lambda
t}u(t)\,dt\,d\lambda\\
&=\int_0^\infty\int_0^{tR}s^{\frac{d}{2}-1}t^{-\frac{d}{2}}e^{-s}u(t)\,ds\,
dt\geq 
\int_1^\infty\int_0^{tR}s^{\frac{d}{2}-1}t^{-\frac{d}{2}}e^{-s}u(t)\,ds\,dt\\
&\geq \left(\int_0^R s^{\frac{d}{2}-1}e^{-s}\,ds\right)\cdot\left(\int_1^\infty
t^{-d/2}u(t)\,dt\right)\,.
\end{align*}
\qed

To handle the case $d \le 2$ in the next proposition and several other places, we will add the following assumption to {\bf (A3)}. 
Note that we do not put the next assumption in Theorem \ref{thm:harnack}.

{\bf (B)}
When $d \le 2$,  we assume that $d+2\delta-2>0$ where $\delta$ is the constant in {\bf (A3)}, and  there are $\sigma'>0$ and  
\begin{equation}\label{e:new22}
\delta'  \in  \left(1-\tfrac{d}{2}, (1+\tfrac{d}{2})
\wedge (2\delta+\tfrac{d-2}{2})\right) 
\end{equation}
 such that
\begin{equation}\label{e:new23}
  \frac{\phi'(\lambda x)}{\phi'(\lambda)}\geq
\sigma'\,x^{-\delta'}\ \text{ for all
}\ x\geq 1\ \text{ and }\ \lambda\geq\lambda_0\,;
\end{equation}

\begin{Prop}\label{prop:pot-2}
If {\bf (A-1)},  {\bf (A-3)} and {\bf (B)} hold and $X$ is transient, then
\[
 g(r)\asymp r^{-d-2}\frac{\phi'(r^{-2})}{\phi(r^{-2})^2},\ r\to 0+\,.
\]
\end{Prop}
\proof
By (\ref{eq:sub-7}) we have 
\[
 g(r)=\int_0^\infty (4\pi
t)^{-d/2}\exp{
\left(-\tfrac{r^2}{4t}\right)}u(t)\,dt\,.
\]
Proposition
\ref{prop:sub-pot} implies that $u(t)\asymp
t^{-2}\frac{\phi'(t^{-1})}{\phi(t^{-1})^2},\ t\to 0+$. 

We are going to use Proposition \ref{prop:gr-1} in Appendix 
with $A=2$, $\eta=u$ and
$\psi=\frac{\phi'}{\phi^2}$. 

In order to use it, we need to check conditions 
(a), (b) and (c)-(ii). The condition (a)
follows from the fact that $\phi'$ and
$\frac{1}{\phi^2}$ are decreasing,
since $\phi$ is a Bernstein function. The condition (b) follows from Lemma
\ref{lem:pot-cond}.

Now we check the condition (c)-(ii) when $d \le 2$; 
Note that by \eqref{e:new22},
$1-\frac{d}{2} < \delta' <2\delta-1 +\frac{d}{2}$ (and $\delta \le 1 <
1+\frac{\delta'}{2}$). Thus 
$0<\delta'+2-2\delta<1+\frac{d}{2}$. Choose $\varepsilon >0$ small so that 
$0<\delta'+2-2\delta+2\varepsilon<1+\frac{d}{2}$, then applying \eqref{e:new23}
and Lemma \ref{lem:sub-50} (ii),  we get 
$$
\frac{\psi(\lambda x)}{\psi(\lambda)}
=\frac{\phi'(\lambda x)}{\phi'(\lambda)}
\frac{\phi(\lambda)^2}{\phi(\lambda x)^2} \ge c_1 x^{-\delta'}c_2 x^{-2+2\delta-2\varepsilon}=c_1c_2x^{-(\delta'+2-2\delta+2\varepsilon)},
$$
Thus \eqref{e:new41} holds. 
\qed

\section{Green function estimates}\label{sec:4}

The purpose of this section is to establish pointwise Green function estimates when $X$ is transient.  
More precisely, 
we are interested in estimate of $G_{B(x_0,r)}(x,y)$ for $x\in
B(x_0,b_1r)$ and $y\in A(x_0,b_2r,r):=\{y\in \R^d\colon
b_2r\leq |y-x_0|<r\}$, for some $b_1, b_2 \in (0,1)$. 
As a starting point we need an estimate of
$G_{B(x_0,r)}(x,y)$
away from the boundary.

In this section,
 we assume that $S=(S_t\colon t\geq 0)$ is
a subordinator with the Laplace exponent $\phi$
 satisfying {\bf (A-1)}--{\bf (A-3)}, {\bf (B)} and assume that 
 $X=(X_t,\P_x)$ is the transient subordinate process
 defined by $X_t=B(S_t)$ where $B=(B_t,\P_x)$ is  a Brownian motion in $\R^d$
independent of $S$.

Recall that, since $X$ is transient,  its
potential $G$ is finite.
\begin{Lem}\label{lem:gr-20}
 There exists $a\in (0,\frac{1}{3})$ and $c_1>0$ such that for any $x_0\in
\R^d$ and 
$r\in (0,1)$
\begin{equation}\label{eq:gr-35}
G_{B(x_0,r)}(x,y)\geq 
c_1\frac{|x-y|^{-d-2}\phi'(|x-y|^{-2})}{\phi(|x-y|^{-2})^2}\ \text{ for all } \
x,y\in
B(x_0,ar)\,.
\end{equation}
In particular, there is a constant $c_2\in (0,1)$ so that
\[
G_{B(x_0,r)}(x,y)\geq  c_2g(|x-y|)\ \text{ for all } \ x,y\in
B(x_0,ar)\,.
\]
\end{Lem}
\proof Let $x,y\in B(x_0,ar)$ with
$a\in (0,1)$ chosen in the course of the proof.
We use (\ref{eq:sub-9}), i.e.
\[
 G_{B(x_0,r)}(x,y)=g(|x-y|)-\E_x[g(|X(\tau_{B(0,r)})-y|)].
\]
Since $|X(\tau_{B(x_0,r)})-y|\geq (1-a)r$ and $|x-y|\leq 2ar$, we get
\[
 |X(\tau_{B(x_0,r)})-y|\geq\tfrac{1-a}{2a}|x-y|.
\]
This together with the fact that $g$ is decreasing yields
\begin{equation}\label{eq:gr-30}
 G_{B(x_0,r)}(x,y)\geq g(|x-y|)-g(\tfrac{1-a}{2a}|x-y|)\,.
\end{equation}
By Proposition \ref{prop:pot-2} there exist constants $0<c_1<c_2$ such that
\begin{equation}\label{eq:gr-31}
 c_1s^{-d+2}\psi(s^{-2})\leq g(s)\leq
c_2s^{-d+2}\psi(s^{-2}),\ s\in (0,1)\,,
\end{equation}
 with
\[
 \psi(\lambda)=\lambda^2\tfrac{\phi'(\lambda)}{\phi(\lambda)^2}, \ \lambda>0\,.
\]
Considering only $a<\frac{1}{3}$ it follows that
$\frac{2a}{1-a}<1$. Combining (\ref{eq:gr-30}), (\ref{eq:gr-31}) 
we arrive at
\begin{align}
& G_{B(x_0,r)}(x,y)\nonumber\\
&\geq c_1|x-y|^{-d+2}\psi(|x-y|^{-2})\left[1-
c_2c_1^{-1}\left(\tfrac{2a}{1-a}\right)^{d-2}\,\frac{\psi\left(\left(\tfrac{2a}{
1-a}
\right)^2|x-y|^ { -2 } \right)}{\psi(|x-y|^{-2})}\right].\label{eq:gr-36}
\end{align}

When $d\ge 3$, 
choose $a<\frac{1}{3}$ small enough so that
$c_2c_1^{-1}\left(\tfrac{2a}{1-a}
\right)^{d-2}\leq \frac{1}{2}$.

When $d \le 2$, using \eqref{e:new22}, first choose $
\varepsilon>0$ small enough so that  $d-2-2\delta'+4\delta-4\varepsilon>0$ then
choose $a<\frac{1}{3}$ small enough so that 
$c_2c_1^{-1}\left(\tfrac{2a}{1-a}
\right)^{d-2-2\delta'+4\delta-4\varepsilon}\leq \frac{1}{2}$.

Then using  the fact that $\lambda  \to \psi(\lambda)$ is increasing (Lemma \ref{lem:gr-10}) when $d \ge 3$, 
and using  \eqref{e:new21}, \eqref{e:new23} and Lemma \ref{lem:sub-50} (ii) when $d \le 2$, 
we get
\begin{align}
&c_2c_1^{-1}\left(\tfrac{2a}{1-a}\right)^{d-2} \frac{\psi\left(\left(\tfrac{2a}{1-a}\right)^2|x-y|^{-2}\right)}{\psi(|x-y|^{-2}
)}
\nonumber
\\
\leq& 
\begin{cases}
 c_2c_1^{-1}\left(\tfrac{2a}{1-a}\right)^{d-2}& \text{when } d \ge 3\nonumber\\
 c_2c_1^{-1}
 \left(\tfrac{2a}{1-a}\right)^{d+2}
 \frac{\phi'\left(\left(\tfrac{2a}{1-a}\right)^2|x-y|^{-2}\right)\phi(|x-y|^{-2}
)^2}{\phi'(|x-y|^{-2}
)\phi\left(\left(\tfrac{2a}{1-a}\right)^2|x-y|^{-2}\right)^2}
 \le 
  \left(\frac{2a}{1-a}\right)^{(d+2)-2\delta'-4+4 \delta -4\varepsilon} 
& \text{when } d \le 2
\end{cases}
\nonumber\\
\le& \frac12. \label{eq:gr-39}
\end{align}
Therefore
   (\ref{eq:gr-36}) and
(\ref{eq:gr-39}) yield
\[
 G_{B(x_0,r)}(x,y)\geq \tfrac{c_1}{2c_2}|x-y|^{-d+2}\psi(|x-y|^{-2}) \text{ for
all }\ x,y\in
B(x_0,ar)\,.
\]
\qed

\begin{Prop}\label{prop:gr-2}
	There exists a constant $c>0$ such that for all $x_0\in \R^d$ and $r\in
(0,1)$
	\[
		\E_x\tau_{B(x_0,r)}\geq \frac{c}{\phi(r^{-2})}\ \text{ for all
}\ x\in B(x_0,\tfrac{ar}{2})\,,
	\]
where $a\in (0,\frac{1}{3})$ 
as
 in Lemma \ref{lem:gr-20}.
\end{Prop}
\proof
	Take $a$ as in Lemma \ref{lem:gr-20} and set $b=\frac{a}{2}$. For any 
$x\in B(x_0,br)$ we have $B(x,br)\subset B(x_0,ar)$ and so it follows from
Lemma \ref{lem:gr-20}  that
	\begin{align*}
		&\E_x\tau_{B(x_0,r)}
		\geq \int_{B(x,br)}G_{B(x_0,r)}(x,y)\,dy\\
		&\geq
c_1\int_{B(x,br)}\frac{|x-y|^{-d-2}\phi'(|x-y|^{-2})}{\phi(|x-y|^{-2})^2}\,dy=\frac{c_2}{\phi(b^{-2}r^{-2})}\geq
\frac{b^2c_2}{\phi(r^{-2})}\,.
	\end{align*}
The last inequality follows Lemma \ref{lem:sub-50}, since $b<1$.
\qed
\begin{Rem}\label{rem:gr-1}
	Note that, by (\ref{eq:sub-9}), for any $x_0\in \R^d$ and $r\in (0,1)$
we
have
	\[
		G_{B(x_0,r)}(x,y)\leq g(|x-y|)\ \text{ for all } \ x,y\in
B(x_0,r)
	\]
	and, consequently, $\E_x\tau_{B(x_0,r)}\leq \frac{c}{\phi(r^{-2})}$ for
any
$x\in B(x_0,r)$.
\end{Rem}

Our approach in obtaining pointwise estimates of  Green function of balls
uses maximum principle for certain operators (in a similar way as in
\cite{BS}).

More precisely, for $r>0$ we define a 
Dynkin-like operator $\sU_r$ by
\[
 (\sU_rf)(x)=\frac{\E_x[f(X(\tau_{B(x,r)}))]-f(x)}{\E_x\tau_{B(x,r)}}
\]
for measurable functions $f\colon\R^d\rightarrow \R$ whenever it is
well-defined.

\begin{Ex}
 Let $x\in \R^d$ and $r>0$. Define 
  \[\eta(z):=\E_z\tau_{B(x,r)},\ z\in\R^d\,.\]
 By the strong Markov property, for any $y\in B(x,r)$ and
$s<r-|y-x|$
\begin{align*}
 \eta(y)=\E_y[\tau_{B(y,s)}+\tau_{B(x,r)}\circ \theta_{\tau_{B(y,s)}}]=\E_y\tau_{B(y,s)}+\E_y\eta(X(\tau_{B(y,s)}))\,.
\end{align*}
Therefore 
\begin{align}\label{e:new2}
(\sU_s \eta)(y)=-1 \quad \text{ for any } y\in B(x,r)  \text{ and }
s<r-|y-x|.
\end{align}
\end{Ex}
\begin{Rem}\label{rem:harm}  Let $h\colon \R^d\rightarrow [0,\infty)$ be a
non-negative
function that is harmonic in bounded open set $D\subset \R^d$. Then for $x\in D$
and $s<\text{dist}(x,\partial D)$ we have $h(x)=\E_x[h(X(\tau_{B(x,s)}))]$. Thus
\[
 (\sU_sh)(x)=0\ \text{ for all }\ x\in D.
\]
\end{Rem}

\begin{Prop}[Maximum principle]\label{prop:max_pr}
 Assume that there exist $x_0\in \R^d$ and $r>0$ such that $(\sU_rf)(x_0)<0$.
Then
\begin{equation}\label{eq:gr-40}
 f(x_0)>\inf_{x\in\R^d}f(x)\,.
\end{equation}
\end{Prop}
\proof
If (\ref{eq:gr-40}) is not true, then $f(x_0)\leq f(x)$ for all $x\in \R^d$.
This implies $(\sU_r f)(x_0)\geq 0$, which is in contradiction with the
assumption. 
\qed

\begin{Prop}\label{prop:green_up}
	There exists a constant $c>0$ such that for
all $r\in (0,1)$ and $x_0\in
\R^d$
	\[
		G_{B(x_0,r)}(x,y)\leq
c\,\frac{r^{-d-2}\phi'(r^{-2})}{\phi(r^{-2})}\,\E_y \tau_{B(x_0,r)}
	\]
	for all $x\in B(x_0,\frac{br}{2})$ and  $y\in A(x_0,br,r)$, where
$b=\frac{a}{2}$ with $a\in (0,\frac{1}{3})$ from Lemma \ref{lem:gr-20}.
\end{Prop}
\proof
	Take $a\in (0,\frac{1}{3})$ as in Proposition \ref{prop:gr-2} (which is
the same one as in  Lemma \ref{lem:gr-20}), 
set $b=\frac{a}{2}$ and let $x\in B(x_0,\frac{br}{2})$ and $y \in A(x_0,br,r)$.
Define functions
\[\eta(z):=\E_z\tau_{B(x_0,r)}\ \text{ and }\ h(z):=G_{B(x_0,r)}(x,z)\]
and choose 
$s<(r-|y|)\wedge\frac{br}{4}$. 
Note that  $h$ is 
harmonic in $B(x_0,r)\setminus\{x\}$ .

	Since $
	h(z)\leq g(\tfrac{br}{8})$  for $ z\in
B(x,\tfrac{br}{8})^c$ 
	and $y \in A(x_0, b r) \subset B(x, br/8)^c$, 
	\eqref{eq:sub-105}, \eqref{e:new1} and Remark \ref{rem:harm} yield
	\begin{align*}
		\sU_s &(h\wedge g(\tfrac{br}{8}))(y)=\sU_s (h\wedge
g(\tfrac{br}{8})-h)(y)
	\\&=\frac{1}{\E_y\tau_{B(y,s)}}\int_{\overline{B(y,s)}^c
} \int_ { B(y ,
s)}G_{B(y , s)}(y,v)j(|z-v|)(h(z)\wedge g(\tfrac{br}{8})-h(z))\,dv\,dz\\
&=\frac{1}{\E_y\tau_{B(y,s)}}\int_{B(x,\frac{br}{8})
} \int_ { B(y ,
s)}G_{B(y , s)}(y,v)j(|z-v|)(h(z)\wedge g(\tfrac{br}{8})-h(z))\,dv\,dz\\
		&\geq
-\frac{1}{\E_y\tau_{B(y,s)}}\int_{B(x,\frac{br}{8})}\int_{B(y,s)}G_{B(y,s)}(y,
v)j(|z-v|)h(z)\,dv\, dz.
	\end{align*}

Note that  $|z-v|\geq  |x-y|-|x-z|-|y-v| \geq \frac{br}{8}$ for $v\in
B(x,\frac{br}{8})$ and $z\in
B(y,s)$ implies $-j(|z-v|) \ge -j(\frac{br}{8})$. Thus
	\begin{align}
		&\sU_s (h\wedge g(\tfrac{br}{8})-h)(y)\nonumber \\
		&\geq
-\frac{j(\tfrac{br}{8})}{\E_y\tau_{B(y,s)}}\left(\int_{B(x,\frac{br}{8})}G_{
B(x_0 ,r)}(x, z)\,
dz\right)\cdot\left(\int_{B(y,s)}G_{B(y,s)}(y,v)\,dv\right)\nonumber\\
		&\geq
-\frac{j(\tfrac{br}{8})}{\E_y\tau_{B(y,s)}}\left(\int_{B(x_0,r)}G_{B(x_0,r)}(x,
z)\, dz\right)\,\E_y\tau_{B(y,s)} \nonumber\\
		&=-j(\tfrac{br}{8})\eta(x)\geq
-c_1\left(\tfrac{b}{8}\right)^{-d-2}\frac{r^{-d-2}\phi'(r^{-2})}{\phi(r^{-2})},
\label{e:new3}
	\end{align}
where in the last inequality we have used Proposition \ref{prop:pot-11}, Remark
\ref{rem:gr-1} and the fact that $\phi'$ is decreasing.

	Similarly, by Proposition \ref{prop:pot-2} and 
Proposition \ref{prop:gr-2}
we see that there is a constant $c_2>0$ such that
	\[
		g(\tfrac{br}{8})\leq
c_2\left(\tfrac{b}{8}\right)^{-d-2}\frac{r^{-d-2}\phi'(r^{-2})}{\phi(r^{-2})}
\eta(z)\ \quad \text{ for all }\ z\in
B(x_0,br)\,.
	\]
		Setting $c_3:=(c_1\vee c_2)\left(\tfrac{b}{8}\right)^{-d-2}+1$
we
obtain
		\[
			c_3\frac{r^{-d-2}\phi'(r^{-2})}{\phi(r^{-2})}\eta(z)
-h(z)\wedge g(\tfrac{br}{8})\geq
c_3\frac{r^{-d-2}\phi'(r^{-2})}{\phi(r^{-2})}\eta(z)-g(\tfrac{br}{8})\geq 0
		\]
		for all $z\in B(x_0,br)$.	
	Therefore, the function
	\[	
u(\cdot):=c_3\frac{r^{-d-2}\phi'(r^{-2})}{\phi(r^{-2})}\eta(\cdot)
-h(\cdot)\wedge
g(\tfrac{br}{8})
	\]
		is nonnegative for  $z\in B(x_0,br)$, vanishes on  $B(x_0,r)^c$
and, by \eqref{e:new2} and \eqref{e:new3},  
	\[
		\sU_s u(y)\leq
-c_3\frac{r^{-d-2}\phi'(r^{-2})}{\phi(r^{-2})}+c_1\left(\tfrac{b}{8}\right)^{
-d-2}\frac{r^{-d-2}\phi'(r^
{ -2 } ) } {
\phi(r^{-2})}<0\ \text{ for } \ y\in A(x_0,br,r)\,.
	\]
 If 
$\inf_{y\in \R^d} u(y)<0$, then by continuity of $u$ on $B(x_0,r)$ there would
exist $y_0\in A(x_0,br,r)$ such
that $u(y_0)= \inf_{y\in \R^d} u(y)$. But then $\sU_su(y_0)\geq 0$,  by the
maximum principle (see Proposition \ref{prop:max_pr}), which is not true.
Therefore $\inf_{y\in \R^d} u(y)\geq 0$.
	
	Finally, since $h\leq g(\tfrac{br}{8})$ on
$A(x_0,br,r)$ it follows that
	\[
		G_{B(x_0,r)}(x,y)\leq c_4
\frac{r^{-d-2}\phi'(r^{-2})}{\phi(r^{-2})} \eta(y)\ \text{ for all }\ y\in
A(x_0,br,r)\,.
	\]
	
\qed

\begin{Prop}\label{prop:green_down}
There exist constants $c>0$ and $b\in (0,1)$ such that for any $r\in
(0,1)$ and $x_0\in \R^d$
\[
	G_{B(x_0,r)}(x,y)\geq
c\,\frac{r^{-d-2}\phi'(r^{-2})}{\phi(r^{-2})}\,\E_y \tau_{B(x_0,r)}
\]
for all $x\in B(x_0,br)$ and $y\in B(x_0,r)$.
\end{Prop}
\proof 
Choose $a\in (0,\frac{1}{3})$ as in Lemma \ref{lem:gr-20}. 
	Then
	\begin{equation}\label{eq:g-2}
		G_{B(x_0,r)}(x,v)\geq
c_1\frac{|x-v|^{-d-2}\phi'(|x-v|^{-2})}{\phi(|x-v|^{-2})^2}\ \text{ for }\
x,v\in
B(x_0,ar)\,.
	\end{equation}
	By Proposition \ref{prop:green_up} we know that there exists a constant
$c_2>0$ so that
	\begin{equation}\label{eq:g-g2}
		G_{B(x_0,r)}(x,v)\leq c_2
\frac{r^{-d-2}\phi'(r^{-2})}{\phi(r^{-2})}\E_v\tau_{B(x_0,r)}\ \text{ for }\
x\in B(x_0,\tfrac{ar}{4}),\ v\in A(x_0,\tfrac{ar}{2},r)\,.
	\end{equation}
	Also, by Remark \ref{rem:gr-1} there is a constant $c_3>0$ such that 
	\begin{equation}\label{eq:g-g3}
		\E_v\tau_{B(x_0,r)}\leq\frac{c_3}{\phi(r^{-2})}\ \text{ for }\
v\in B(x_0,r)\,.
	\end{equation}

	 We take
	\[
		b\leq
\min\left\{\tfrac{1}{4}\left(\tfrac{c_1}{2c_2c_3}\right)^{1/d},\tfrac{a}{8}
\right\}
	\]
and fix it.
	Note that $c_2c_3\leq \frac{c_1}{2}(4b)^{-d}$, i.e. $b\leq
\frac{1}{4}\left(\frac{c_1}{2c_2c_3}\right)^{1/d}$. Thus by Lemma
\ref{lem:gr-10}
	\[
		c_2c_3\frac{r^{-d-2}\phi'(r^{-2})}{\phi(r^{-2})^2}\leq
\tfrac{c_1}{2}\frac{(4br)^{-d-2}\phi'((4b)^{-2}r^{-2})}{\phi((4b)^{-2}r^{-2})^2}
\,.
	\] 
	Now, by (\ref{eq:g-2}) and  (\ref{eq:g-g3}), for all $x\in B(x_0,br)$
and  $v\in B(x,br)$
	\begin{equation}\label{eq:g-g4}
		c_2
\frac{r^{-d-2}\phi'(r^{-2})}{\phi(r^{-2})}\E_v\tau_{B(x_0,r)}\leq
\tfrac{1}{2}G_{B(x_0,r)}(x,v)\,.
	\end{equation}

	For the rest of the proof, we fix $x \in B(x_0,br)$ and 
	define a function 
	\[
		h(v)=G_{B(x_0,r)}(x,
		v)\wedge \left(
c_2\frac{r^{-d-2}\phi'(r^{-2})}{\phi(r^{-2})}\E_v\tau_{B(x_0,r)}\right)\,.
	\]
	
	Let $y\in A(x_0,\frac{ar}{2},r)$ and take $s<(r-|y|)\wedge
\frac{br}{8}$\,. Note that, by (\ref{eq:g-g4}),
	\[
		h(v)\leq \tfrac{1}{2}G_{B(x_0,r)}(x,v)\ \text{ for }\ v\in
B(x,br).
	\]
	Therefore,   (\ref{eq:sub-105}) and Remark \ref{rem:harm} yield
	\begin{align*}
		(\sU_s &h)(y)=\sU_s
\left(h-G_{B(x_0,r)}(x,\cdot)\right)(y)	\\
		&=\frac{1}{\E_y
\tau_{B(y,s)}}\int_{\overline{B(y,s)}^c}\int_{B(y,s)}G_{B(y,s)}(y,
v)j(|z-v|) \left(h(z)-G_{B(x_0,r)}(x,z)\right)\,dv\,dz\\
	&\leq
-\frac{1}{2\E_y
\tau_{B(y,s)}}\int_{B(x,br)}\int_{B(y,s)}G_{B(y,s)}(y,
v)j(|z-v|)G_{B(x_0,r)}(x,z)\, dv\,dz\,.
	\end{align*}
	Note that in the second equality we have used that
$h(y)=G_{B(x_0,r)}(x,y)$, which follows from (\ref{eq:g-g2})\,.
	
Since $|z-v|\leq 2r$ we obtain
	\begin{align*}
		(\sU_s h)(y)
		&\leq-\frac{j(2r)}{2\E_y\tau_{B(y,s)}}
		\left(\int_{B(x,br)}G_{B(x_0,r)}(x,z)\,dz\right)\,\E_y
\tau_{B(y,s)}
	\end{align*}
By Proposition \ref{prop:pot-11}, (\ref{eq:g-2}) and
the fact that $\lambda\mapsto\frac{\lambda}{\phi(\lambda)}$ 
is increasing
(by Lemma \ref{lem:sub-50} or using the fact that $\phi$ is a complete Bernstein
function)
 and $\phi'$ decreasing we finally arrive
at
	\begin{align*}
		(\sU_s h)(y)\leq
-c_4\frac{r^{-d-2}\phi'(r^{-2})}{\phi\left(\left(br\right)^{-2}
\right) }\leq -c_5 \frac{r^{-d-2}\phi'(r^{-2})}{\phi(r^{-2})}\,.
	\end{align*}

	Define $u=h-\kappa\eta$,
where $\eta(v)=\E_v\tau_{B(x_0,r)}$ and 
	\[
		\kappa=\min\left\{
\frac{c_5}{2},\frac{c_1}{2c_3},\frac{c_2}{2}\right\}\frac{r^{-d-2}\phi'(r^{-2})}
{\phi(r^{-2})}\,.
	\]
	For $y\in A(x_0,\frac{ar}{2},r)$ 
	we have by \eqref{e:new2},  
	\[
		(\sU_s u)(y)\leq 
-c_5\frac{r^{-d-2}\phi'(r^{-2})}{\phi(r^{-2})}+\kappa\leq
-\tfrac{c_5}{2}\frac{r^{-d-2}\phi'(r^{-2})}{\phi(r^{-2})}<0\,.
	\]
	On the other hand, by
	\eqref{eq:g-2}--\eqref{eq:g-g3}, 	for all $v\in B(x_0,\frac{ar}{2})$,  
	\begin{align*}
		u(v)&\geq \left(\tfrac{c_1}{c_3}\wedge
c_2\right)\frac{r^{-d-2}\phi'(r^{-2})}{\phi(r^{-2})}\E_v\tau_{B(x_0,r)}-\kappa\,
\E_v\tau_{B(x_0,r)}\\&\geq
\left(\tfrac{c_1}{2c_3}\wedge
\tfrac{c_2}{2}\right)\frac{r^{-d-2}\phi'(r^{-2})}{\phi(r^{-2})}\E_v\tau_{B(x_0,
r)}\geq 0\,.
	\end{align*}

	Similarly as in Proposition \ref{prop:green_up}, by the maximum
principle it follows that
	\[
		u(y)\geq 0\ \text{ for all }\ y\in B(x_0,r)\,.	
	\]	
\qed

Combining Propositions \ref{prop:green_up} and  \ref{prop:green_down}
we obtain an  important estimate for the Green function.
\begin{Cor}\label{cor:green}
 There exist constants $c_1,c_2>0$ and $b_1,b_2\in (0,\frac{1}{2})$, $2b_1<b_2$
such that for all
$x_0\in \R^d$ and $r\in (0,1)$
\[
	c_1 \frac{r^{-d-2}\phi'(r^{-2})}{\phi(r^{-2})}\,\E_y\tau_{B(x_0,r)}\leq
G_{B(x_0,r)}(x,y)\leq c_2
\frac{r^{-d-2}\phi'(r^{-2})}{\phi(r^{-2})}\,\E_y\tau_{B(x_0,r)}\ 
\]
for all $x\in B(x_0,b_1r)$ and $y\in A(x_0,b_2r,r)$.
\qed
\end{Cor}

\section{Poisson kernel and Harnack inequality}\label{sec:5}

The first goal of this section is to estimate 
Poisson kernel of $X$ in a ball
 given
by
\begin{equation}\label{eq:pkhi-10}
	K_{B(x_0,r)}(x,z)=\int_{B(x_0,r)}G_{B(x_0,r)}(x,y)j(|z-y|)\,dy\,,
\end{equation}
where $x\in 
B(x_0,br)$
 and $z\not\in B(x_0,r)$.
\begin{Prop}\label{prop:poisson}

Suppose that $\phi$
 satisfying {\bf (A-1)}--{\bf (A-3)}, {\bf (B)} and that 
 the corresponding subordinate Browninan motion $X=(X_t,\P_x)$ is  transient.
Then there exist constants $c>0$ and $b\in (0,1)$ such that for all
$x_0\in \R^d$ and $r\in (0,1)$ 
	\[
		K_{B(x_0,r)}(x_1,z)\leq c K_{B(x_0,r)}(x_2,z)
	\]
	for all $x_1,x_2\in B(x_0,br)$ and $z\in B(x_0,r)^c$.
\end{Prop}
\proof
  Take $b_1,b_2\in (0,\frac{1}{2})$ as in Corollary \ref{cor:green}, and let
$x_0\in \R^d$, $x_1,x_2\in
B(x_0,b_1r)$ and $z\in B(x_0,r)^c$.

	We split the integral in (\ref{eq:pkhi-10}) in two
parts
	\begin{align*}	
K_{B(x_0,r)}(x_1,z)=&\int_{B(x_0,b_2r)}G_{B(x_0,r)}(
x_1,y)j(|z-y|)\,dy\,
\\&+ \int_{A(x_0,b_2r,r)}G_{B(x_0,r)}(
x_1,y)j(|z-y|)\,dy\,=:\,I_1+I_2\,.
	\end{align*}

To estimate $I_2$ we use Corollary \ref{cor:green} to 
get that for $y\in A(x_0,b_2 r,r)$
	\[
		c_1
\frac{r^{-d-2}\phi'(r^{-2})}{\phi(r^{-2})}\,\E_y\tau_{B(x_0,r)}\leq
G_{B(x_0,r)}(
x_i,y)\leq c_2 \frac{r^{-d-2}\phi'(r^{-2})}{\phi(r^{-2})}\,\E_y\tau_{B(x_0,r)},
\quad i=1,2.
	\]
	Therefore
	\begin{align*}	
I_2&=\int_{A(x_0,b_2r,r)}G_{B(x_0,r)}(x_1,y)j(|z-y|)\,dy\\
			&\leq
c_1\frac{r^{-d-2}\phi'(r^{-2})}{\phi(r^{-2})}\int_{A(x_0,b_2r,r)}\E_y\,
\tau_ { B(x_0,r)}j(|z-y|)\,dy\\
			&\leq \frac{c_1}{c_2}
\int_{A(x_0,b_2r,r)}G_{B(x_0,r)}(x_2,y)j(|z-y|)\,dy\,\leq \, \frac{c_1}{c_2}K_{B(x_0,r)}(x_2,z)\,.
  \end{align*}	

To handle $I_1$ we consider two cases. If $z\in
A(x_0,r,2)$, then
\[
(1-b_2)|z-x_0|\leq |z-y|\leq
(1+b_2)|z-x_0| \ \text{ for all }\ y\in
B(x_0,b_2r)\,.
\]
Since $1-b_2\geq \frac{1}{2}$ and $1+b_2\leq 2$ we obtain
\begin{equation}\label{eq:hi-1}
 j\left(2|z-x_0|\right)\leq j(|z-y|)\leq j\left(\tfrac{1}{2}|z-x_0|\right)\,.
\end{equation}
Using  Proposition \ref{prop:pot-11} and the fact that $\phi'$ is decreasing
we see that 
\begin{equation}\label{eq:hi-2}
j\left(\tfrac{1}{2}|z-x_0|\right)\leq c_3j\left(2|z-x_0|\right)\,.
\end{equation}

Lemma \ref{lem:gr-20}, Remark \ref{rem:gr-1}, (\ref{eq:hi-1}) and
(\ref{eq:hi-2}) yield
\begin{align*}
 I_1&\leq
j\left(\tfrac{1}{2}|z-x_0|\right)\int_{B(x_0,b_2r)}G_{B(x_0,r)}(x_1,y)\,
dy\\
&\leq c_3j\left(2|z-x_0|\right)\frac{c_4}{\phi((b_2r)^{-2})}\\
&\leq c_5 \int_{B(x_0,b_2r)}G_{B(x_0,r)}(x_2,y)j(|z-y|)\,dy\,\leq \, c_5 K_{B(x_0,r)}(x_2,
z)\,. 
\end{align*}

When $z\in B(x_0,2)^c$ we use
$
 |z-x_0|-br\leq |z-y|\leq |z-x_0|+br$ 
  for all $\ y\in B(x_0,br)$.
Since $b_2\in (0,\frac{1}{2})$ and $r\in (0,1)$ we have
\begin{equation}\label{eq:hi-3}
 j(|z-x_0|+\tfrac{1}{2})\leq j(|z-y|)\leq j(|z-x_0|-\tfrac{1}{2})\,.
\end{equation}
By (\ref{eq:sub-11}) we have
\begin{equation}\label{eq:hi-4}
 j(|z-x_0|-\tfrac{1}{2})\leq c_6 j(|z-x_0|+1)\leq
j(|z-x_0|+\tfrac{1}{2})\,.
\end{equation}

Similar to the previous case, by (\ref{eq:sub-11}), Lemma \ref{lem:gr-20},
(\ref{eq:hi-3}) and
(\ref{eq:hi-4})  we have
$
 I_2\leq c_7K_{B(x_0,r)}(x_2,z)$. 
Therefore,
$
 I_2\leq (c_5\vee c_7)K_{B(x_0,r)}(x_2,z)
$
which implies that 
\[
 K_{B(x_0,r)}(x_1,z)\leq
\max\{c_5,c_7,\tfrac{c_1}{c_2}\}K_{B(x_0,r)}(x_2,z)\,.
\]
\qed

Now we are ready to prove our main result.

\proof[Proof of Theorem \ref{thm:harnack}]
	Suppose $d \ge 3$. Then $X$ is always transient. Take $b>0$ as in Proposition \ref{prop:poisson} and set
$a=\frac{b}{4}$. 
Suppose that $h\colon \R^d\rightarrow
[0,\infty)$ is harmonic in $B(x_0,r)$.
Using representation 
$$
h(x)=\E_x[h(X_{\tau_{B(x_0,2ar)}})]=\int_{\overline{B(x_0,2ar)}^c}K_{B(x_0,2ar)}
(x_1,
z)h(z)\, dz $$
 and Proposition
\ref{prop:poisson} we have
	\begin{align*}	
h(x_1)&=\int_{\overline{B(x_0,2ar)}^c}K_{B(x_0,2ar)}(x_1,
z)h(z)\, dz\\&\leq c
\int_{\overline{B(x_0,2ar)}^c}K_{B(x_0,2ar)}(x_2,z)h(z)\,
dz=c\,h(x_2).
	\end{align*}
	Then, using a standard Harnack chain argument, we prove the theorem for
$d\ge 3$.
	
	To handle the lower dimensional case, we use the following notation:
for $x=(x^1,\ldots,x^{d-1},x^d)\in \R^d$ we set
$\tilde{x}=(x^1,\ldots,x^{d-1})$.
Let $X=((\widetilde{X}_t,X^d_t), \P_{(\widetilde{x}, x^d)})$
 be a $d$-dimensional subordinate Brownian motion
with the characteristic exponent 
\[
 \Phi(\xi)=\phi(|\xi|^2),\ \xi\in\R^d\,.
\]
 By checking the
characteristic functions, it follows that, for every $x^d \in \R$, 
$\widetilde{X}=(\widetilde{X}_t,\P_{\widetilde{x}})$ is a $(d-1)$-dimensional
subordinate Brownian motion with characteristic exponent
\[
 \widetilde{\Phi}(\tilde{\xi})=\phi(|\tilde{\xi}|^2),\
\tilde{\xi}\in \R^{d-1}\,.
\]

Suppose the theorem is true for 
for some $d\geq 2$.
 Let $h\colon \R^{d-1} \rightarrow
[0,\infty)$ be a function that is harmonic in $B(
\wt {x}_0,r)$. 

Since \[\tau_{\B(\wt {x}_0,s) \times \R}=\inf\{t>0 : \wt
{X}_t \notin \B(\wt {x}_0,s)\} ,\]  the strong Markov property   implies that
the function $f \colon \R^{d} \rightarrow
[0,\infty)$ defined by $f(\wt x, x^d)=h(\wt x)$ is harmonic in $B(
\wt {x}_0,r) \times \R$. 

In particular,  $f$ is harmonic in $B((\wt {x}_0, 0),r)$.
By applying the result to  $f$, we see that there exists a constant $c>0$  such
that for
all $\wt {x}_0 \in \R^{d-1}$ and $r\in (0,1)$
	\[
		h(\wt {x}_1)=f((\wt {x}_1, 0))\leq c\, f((\wt {x}_2, 0))=
c\ h(\wt {x}_2) \ \text{ for all }\ \wt{x}_1,\wt{x}_2\in
B(\wt {x}_0,\tfrac{r}{2}).
	\]

Applying this argument first to $d=3$ and then to $d=2$ we finish the proof of
the theorem.
\qed

Since
$K_{B(x_0,r)}(x,\cdot)$  is continuous on $\overline{B(x_0,r)}^{\, c}$ for every
$x\in B(x_0,r)$, Theorem \ref{thm:harnack} implies Proposition
\ref{prop:poisson} without the conditions {\bf (B)} and $X$ being transient.

\begin{Cor}\label{c:poisson}

Suppose that $\phi$
 satisfying {\bf (A-1)}--{\bf (A-3)}. 
Then for every $b\in (0,1)$,  there exists a constant $c=c(b)>0$  such that for
all
$x_0\in \R^d$ and $r\in (0,1)$ 
	\[
		K_{B(x_0,r)}(x_1,z)\leq c K_{B(x_0,r)}(x_2,z)
	\]
	for all $x_1,x_2\in B(x_0,br)$ and $z\in B(x_0,r)^c$.
\end{Cor}

We
we omit the proof since it is the same as the proof of 
Proposition 1.4.11 in \cite{KSV3}.

\appendix
\section{Asymptotical properties}\label{app:a}
In this section we always assume that $f\colon (0,\infty)\rightarrow (0,\infty)$
is  a differentiable function
satisfying 
	\begin{equation}\label{eq:cond_main}
	|f(\lambda+\varepsilon)-f(\lambda)|=\int_0^\infty \left(e^{-\lambda
t}-e^{-(\lambda+\varepsilon) t}\right)\,\nu(t)\,dt\,,
\end{equation}
for all $\lambda>0$, $\varepsilon\in (0,1)$
and a decreasing function   $\nu\colon (0,\infty)\rightarrow (0,\infty)$.

\begin{Lem}\label{lem:1} 
For all $t>0$, 
	\[
		\nu(t)\leq (1-2e^{-1})^{-1}\,t^{-2}|f'(t^{-1})|\,.
	\]
\end{Lem}
\proof
	Let $\varepsilon\in (0,1)$. Then
	\begin{align*}
		|f(\lambda+\varepsilon)-f(\lambda)|&=\int_0^\infty
\left(e^{-\lambda t}-e^{-\lambda t-\varepsilon t}\right)\nu(t)\,dt\\
		&=\lambda^{-1}\int_0^\infty
e^{-z}\left(1-e^{-\varepsilon\lambda^{-1}
z}\right)\nu(\lambda^{-1}z)\,dz\,.
	\end{align*}
	
	Since $\nu$ is decreasing, for any $r>0$ we conclude
	\begin{align*}
		|f(\lambda+\varepsilon)-f(\lambda)|&\geq \lambda^{-1}\int_0^r
e^{-z}\left(1-e^{-\varepsilon
\lambda^{-1}z}\right)\nu(\lambda^{-1}z)\,dz\\
		&\geq \lambda^{-1}\nu(\lambda^{-1}r)\int_0^r
e^{-z}\left(1-e^{-\varepsilon\lambda^{-1} z}\right)\,dz.
	\end{align*}
	Therefore
	\begin{align}\label{eq:1}	
\left|\frac{f(\lambda+\varepsilon)-f(\lambda)}{\varepsilon}\right|& \geq
\lambda^{-2}\nu(\lambda^{-1}r)\int_0^r
z\,e^{-z}\,\frac{1-e^{-\varepsilon \lambda^{-1}z}}{\varepsilon
\lambda^{-1}z}\,dz.
	\end{align}
	By
	the Fatou's lemma and (\ref{eq:1}) we obtain
	\begin{align*}
		|f'(\lambda)|&=\lim\limits_{\varepsilon \to
0+}\left|\frac{f(\lambda+\varepsilon)-f(\lambda)}{\varepsilon}\right|
		\geq \lambda^{-2}\nu(\lambda^{-1}r)\int_0^r ze^{-z}\,dz\\
		&=
\lambda^{-2}\nu(\lambda^{-1}r)\left(1-e^{-r}(r+1)\right)\,.
	\end{align*}
	
	In particular, for $r=1$ we deduce
	\[
		\nu(t)\leq \left(1-2e^{-1}\right)^{-1}t^{-2}|f'(t^{-1})|,\
t>0\,.
	\]
\qed

\begin{Lem}\label{lem:2}
	Assume that $|f'|$ is decreasing and there exist $c_1>0$,
$\lambda_0>0$ and $\delta>0$ such that
	\begin{equation}\label{eq:cond}
		\left|\frac{f'(\lambda x)}{f'(\lambda)}\right|\leq
c_1x^{-\delta}\ \text{ for all }\ \lambda\geq \lambda_0\ \text{ and }\ x\geq
1\,.
	\end{equation}
	Then there is a constant $c_2=c_2(c_1,\lambda_0,\delta)>0$ such
that 
	\[
		\nu(t)\geq c_2t^{-2}|f'(t^{-1})|\ \text{  for any }\ t\leq
1/\lambda_0\,.
	\]
\end{Lem}
\proof
Let $\varepsilon \in (0,1)$. For $r\in (0,1]$ we have 
	\begin{align}\nonumber
		|f(\lambda+\varepsilon)-f(\lambda)|
		&=\lambda^{-1}\int_0^\infty e^{-z}\left(1-e^{-\varepsilon
\lambda^{-1}z}\right)\nu(\lambda^{-1}z)\,dz\\
		&=I_1(\varepsilon)+I_2(\varepsilon)\,,\label{eq:t1}
	\end{align}
	where
	\begin{align*}
		I_1(\varepsilon)&=\lambda^{-1}\int_0^r
e^{-z}\left(1-e^{-\varepsilon \lambda^{-1}z}\right)\nu(\lambda^{-1}z)\,dz\\
		I_2(\varepsilon)&=\lambda^{-1}\int_r^\infty
e^{-z}\left(1-e^{-\varepsilon \lambda^{-1}z}\right)\nu(\lambda^{-1}z)\,dz\,.
	\end{align*}
		
	Since $\nu$ is decreasing, 
	\begin{align*}
		\frac{I_2(\varepsilon)}{\varepsilon}\leq
\lambda^{-2}\nu(\lambda^{-1}r)\int_r^\infty
ze^{-z}\,\frac{1-e^{-\varepsilon \lambda^{-1}z}}{\varepsilon
\lambda^{-1}z}\,dz\,,
	\end{align*}
	 and so by the dominated convergence theorem we deduce
	\begin{equation}\label{eq:3}
		\limsup_{\varepsilon \to
0+}\frac{I_2(\varepsilon)}{\varepsilon}\leq
\lambda^{-2}\nu(\lambda^{-1}r)\int_r^\infty
ze^{-z}\,dz=(r+1)e^{-r}\lambda^{-2}\nu(\lambda^{-1}r)\,.
	\end{equation}
	
	On the other hand, by Lemma \ref{lem:1} and (\ref{eq:cond}) we have
	\begin{align*}
		\frac{I_1(\varepsilon)}{\varepsilon}&\leq
\frac{\lambda^{-2}}{1-2e^{-1}}\int_0^r ze^{-z}\,\frac{1-e^{-\varepsilon
\lambda^{-1}z}}{\varepsilon
\lambda^{-1}z}\,\frac{|f'(\lambda z^{-1})|}{\lambda^{-2}z^2}\,dz\\
		&\leq \frac{c_1}{1-2e^{-1}}\,|f'(\lambda)|\int_0^r
e^{-z}\,\frac{1-e^{-\varepsilon \lambda^{-1}z}}{\varepsilon
\lambda^{-1}z}z^{\delta-1}\,dz\,.
	\end{align*}
	The dominated convergence implies
	\begin{equation}\label{eq:4}
	\limsup_{\varepsilon \to 0+} \frac{I_1(\varepsilon)}{\varepsilon}\leq
\frac{c_1}{1-2e^{-1}}\,|f'(\lambda)|\int_0^r e^{-z}z^{\delta -1}\,dz\,\ \text{
for any }\ \lambda\geq \lambda_0.
	\end{equation}
	
	Combining (\ref{eq:t1}), (\ref{eq:3}) and (\ref{eq:4}) we deduce
	\[
		|f'(\lambda)|\leq  \frac{c_1}{1-2e^{-1}}\,|f'(\lambda)|\int_0^r
e^{-z}z^{\delta -1}\,dz+(r+1)e^{-r}\lambda^{-2}\nu(\lambda^{-1}r)
	\]
	for all $\lambda\geq \lambda_0$.
	
	Choosing $r_0\in (0,1]$ so that
	\[
	\frac{c_1}{1-2e^{-1}}\,\int_0^{r_0} e^{-z}z^{\delta -1}\,dz\leq
\frac{1}{2}\,.
	\]
	we have
	\[
		\nu(\lambda^{-1}r_0)\geq
\frac{e^{r_0}}{2(r_0+1)}\lambda^2|f'(\lambda)|\ \text{ for all }\ \lambda\geq
\lambda_0.
	\]
	Since $|f'|$ is decreasing, we see that
	\begin{align*}
		\nu(t)&\geq 
\frac{e^{r_0}}{2(r_0+1)}\frac{|f'(r_0/t)|}{(t/r_0)^2}\\
			&\geq
\frac{r_0^2e^{r_0}}{2(r_0+1)}t^{-2}|f'(t^{-1})|\ \text{ for all }\ t\leq
r_0\lambda_0^{-1}.
	\end{align*}
	On the other hand, for 	$r_0\lambda_0^{-1}\leq t\leq \lambda_0^{-1}$ we
have
\[
  \nu(t)\geq \nu(\lambda_0^{-1})\geq
t^{-2}|f'(t^{-1})|\frac{(r_0/\lambda_0)^2}{|f'(\lambda_0)|},
\]
since $\nu$ and $|f'|$ are decreasing.

	Setting \[c_2=\frac{r_0^2e^{r_0}}{2(r_0+1)}\wedge\frac{
\nu(\lambda_0^{-1})
\lambda_0^{-2}r_0^2}{|f'(\lambda_0)|}\] we get
	\[
		\nu(t)\geq c_2 t^{-2}|f'(t^{-1})|\ \text{ for all }\ t\leq
\lambda_0^{-1}\,.
	\]	
\qed

\begin{Prop}\label{prop:gr-1}
Let $A>0$ and $\eta\colon (0,\infty)\rightarrow (0,\infty)$ be a decreasing
function
satisfyng the following conditions:
\begin{itemize}
 \item[(a)] there exists a decreasing function $\psi\colon (0,\infty)\rightarrow
(0,\infty)$ such that $\lambda\mapsto \lambda^2\psi(\lambda)$ is increasing and
satisfies \[\eta(t)\asymp 
t^{-A} \psi(t^{-1}),\ t\to 0 +\, ;\]
 \item[(b)] $\int_1^\infty t^{-d/2}\eta(t)\,dt<\infty$
\item[(c)] 
either 
 (i)
$A> 3-\frac{d}{2}$ 
or 
(ii) 
$A> 3-\frac{d}{2}$  when $d\geq 3 $ and in the case $d \le 2$ there exist
$\delta>0$ and $c>0$ such that $A-\delta>1-\frac{d}{2}$ and 
\begin{align}\label{e:new41}
 \frac{\psi(\lambda x)}{\psi(\lambda)}\geq cx^{-\delta}\ \text{ for all }\
x\geq 1\ \text{ and }\ \lambda\geq 1\,.
\end{align}
\end{itemize}
If 
\[
I(r)=\int_0^\infty (4\pi
t)^{-d/2}\exp{\left(-\frac{r^2}{4t}\right)}\eta(t)\,dt
\]
exists for $r\in (0,1)$ small enough, then
\[
 I(r)\asymp 
 r^{-d-2A+2}
 \psi(r^{-2}),\ r\to 0+\,.
\]
\end{Prop}
\proof
Write for $r>0$ 
\begin{align}\nonumber
 I(r)&=\int_0^{r^2}(4\pi t)^{-d/2}\exp{\left(-\tfrac{r^2}{4t}\right)}
\eta(t)\,dt+\int_{r^2}^\infty (4\pi
t)^{-d/2}\exp{\left(-\tfrac{r^2}{4t}\right)}\eta(t)\,dt\\&=I_1(r)+I_2(r)\,
.\label{eq:gr-1}
\end{align}
By condition (a), 
\begin{align}
 I_1(r)&\leq c_1\int_0^{r^2}(4\pi
t)^{-d/2}\exp{\left(-\tfrac{r^2}{4t}\right)}
t^{-A}
\psi(t^{-1})\,dt\nonumber\\&\leq
c_2\psi(r^{-2})\int_0^{r^2}t^{-\frac{d}{2}-
A}\exp{\left(-\tfrac{
r^2}{4t}\right)}\,dt\nonumber\\
&=c_3r^{-d-2A+2}\psi(r^{-2})\int_{\frac{1}{4}}^\infty
t^{A-2+\frac{d}{2}}e^{-t}\,dt\,\label{eq:gr-2}\,.
\end{align}
Similarly,  
\begin{align}
 I_2(r)&\leq c_1\int_{r^2}^1(4\pi
t)^{-d/2}\exp{\left(-\tfrac{r^2}{4t}\right)}
t^{-
A}\psi(t^{-1})\,dt+\int_1^\infty(4\pi
t)^{-d/2}\exp{\left(-\tfrac{r^2}{4t}\right)}
\eta(t)\,dt\nonumber\\
&\leq c_1\int_{r^2}^1(4\pi
t)^{-d/2}
t^{-
A}\psi(t^{-1})\,dt+\int_1^\infty(4\pi
t)^{-d/2}\eta(t)\,dt\nonumber\,.
\end{align}
The following inequality holds
 \begin{equation}\label{eq:tmp-aymp-13}
  \int_{r^2}^1(4\pi
t)^{-d/2}
t^{-A}\psi(t^{-1})\,dt \leq c_4 r^{-d-2A+2}\psi(r^{-2})\,,
 \end{equation}
since 
\begin{itemize} 
\item[(1)]
if condition (c)-(i) holds, then by conditions (a) and (c)-(i)
\begin{align*}
 \int_{r^2}^1(4\pi
t)^{-d/2}
t^{-
A}\psi(t^{-1})\,dt&\leq r^{-4}\psi(r^{-2})\int_{r^2}^1(4\pi
t)^{-d/2}t^{2-A}\,dt\\&\leq c_5 r^{-d-2A+2}\psi(r^{-2})\,;
\end{align*}
\item[(2)]
if condition (c)-(ii) holds and  $d \le 2$, \eqref{e:new41} implies
\begin{align*}
 \int_{r^2}^1(4\pi
t)^{-d/2}
t^{-
A}\psi(t^{-1})\,dt&\leq
\psi(r^{-2})r^{-2\delta}
\int_{r^2}^1(4\pi t)^{-\frac{d}{2}}t^{
-A+\delta}\,dt\\
&\leq c_6 r^{-d-2A+2}\psi(r^{-2})\,.
\end{align*}
\end{itemize}
In particular,  (\ref{eq:tmp-aymp-13}) implies
\[
 r^{-d-2A+2}\psi(r^{-2})\geq c_7>0\ \text{ for all }\ r\in (0,1)\,
\]
and thus
\begin{align}
I_2(r)&
\leq
c_6r^{-d-2A+2}
\psi(r^{-2})+c_8
\leq c_{9}
r^{-d-2A+2}
\psi(r^{-2})\label{eq:gr-3}\,.
\end{align}
Combining (\ref{eq:gr-1}), (\ref{eq:gr-2}) and (\ref{eq:gr-3}) we get the upper
bound
\[
 I(r)\leq c_7 
 r^{-d+2-2A}\psi(r^{-2})\ \text{ for all } \ r\in (0,1)\,.
\]
To get the lower bound we estimate $I(r)$ from below by $I_1(r)$ and use
(a) to get
\begin{align*}
 j(r)&\geq I_1(r)\geq c_8 \int_0^{r^2}(4\pi
t)^{-d/2}\exp{\left(-\tfrac{r^2}{4t}\right)}
t^{-
A}\psi(t^{-1})\,dt\nonumber\\
&\geq c_8 r^{-4}\psi(r^{-2})\int_0^{r^2}
(4\pi t)^{-d/2} t^{2-A}
\exp{\left(-\tfrac{r^2}{4t}\right)}\,dt\nonumber\\
&=c_9r^{-d-2A+2}
\psi(r^{-2})\int_{\frac{1}{4}}^\infty
s^{-\frac{d}{2}+A-4}e^{-s}\,ds\\&=c_{10}
r^{-d-2A+2}
\psi(r^{-2})\ \text{ for all }\
r\in  (0,1).
\end{align*}
\nopagebreak[4]
\qed

\providecommand{\bysame}{\leavevmode\hbox to3em{\hrulefill}\thinspace}
\providecommand{\MR}{\relax\ifhmode\unskip\space\fi MR }
\providecommand{\MRhref}[2]{%
  \href{http://www.ams.org/mathscinet-getitem?mr=#1}{#2}
}
\providecommand{\href}[2]{#2}

\end{document}